# Some results on the Gittins index for a normal reward process

Yi-Ching Yao[1]

*Academia Sinica*

**Abstract:** We consider the Gittins index for a normal distribution with unknown mean $\theta$ and known variance where $\theta$ has a normal prior. In addition to presenting some monotonicity properties of the Gittins index, we derive an approximation to the Gittins index by embedding the (discrete-time) normal setting into the continuous-time Wiener process setting in which the Gittins index is determined by the stopping boundary for an optimal stopping problem. By an application of Chernoff's continuity correction in optimal stopping, the approximation includes a correction term which accounts for the difference between the discrete and continuous-time stopping boundaries. Numerical results are also given to assess the performance of this simple approximation.

## 1. Introduction

The classical multi-armed bandit problem is concerned with sequential design of adaptive sampling from $k$ statistical populations with distribution functions $F_{\theta_i}$, $i = 1, \ldots, k$ ($k \geq 2$) where $\theta_i$ denotes the unknown parameter of the $i$th population. Specifically, the objective is to sample $Y_1, Y_2, \ldots$ sequentially from the $k$ populations so as to maximize the expected total discounted reward

$$E_\pi E_{\theta_1,\ldots,\theta_k}\left(\sum_{j=1}^\infty \gamma_j Y_j\right) = \int E_{\theta_1,\ldots,\theta_k}\left(\sum_{j=1}^\infty \gamma_j Y_j\right) d\pi(\theta_1, \ldots, \theta_k),$$

where $\pi$ is the prior distribution of $(\theta_1, \ldots, \theta_k)$ and $\{\gamma_j\}$ is a (deterministic) discount sequence. The two most important types of discount sequence are uniform discounting with finite horizon $N > 0$ (i.e. $\gamma_j = 1$ for $j \leq N$ and $\gamma_j = 0$ for $j > N$) and geometric discounting with discount factor $0 < \beta < 1$ (i.e. $\gamma_j = \beta^{j-1}, j = 1, 2, \ldots$). While in general the optimal allocation rule can only be characterized via the dynamic programming equations which admit no general closed-form solutions, Gittins and Jones [13] showed that under geometric discounting, when the prior distribution is a product measure $d\pi(\theta_1, \ldots, \theta_k) = d\pi_1(\theta_1) \times \cdots \times d\pi_k(\theta_k)$, the optimal allocation rule is to sample at each stage from the population with the greatest (current) Gittins index. See also [11] and [21].

For a population with distribution function $F_\theta$ and (current) prior distribution $\pi(\theta)$ of the unknown parameter $\theta$, the Gittins index is defined by

$$(1) \qquad \lambda(\pi, \beta) = \sup_{\xi \geq 1}\left[E_\pi E_\theta\left(\sum_{n=1}^\xi \beta^{n-1} X_n\right) \bigg/ E_\pi E_\theta\left(\sum_{n=1}^\xi \beta^{n-1}\right)\right]$$

[1]Institute of Statistical Science, Academia Sinica, Taipei, Taiwan, e-mail: yao@stat.sinica.edu.tw

where the supremum is taken over all (integer-valued) stopping times $\xi \geq 1$ and $X_1, X_2, \ldots$ are (conditionally) iid with common distribution function $F_\theta$ (given $\theta$). Equivalently, $\lambda(\pi, \beta)$ is the infimum of the set of solutions $\lambda$ of the equation

$$(2) \qquad \frac{\lambda}{1-\beta} = \sup_{\xi \geq 0} E_\pi E_\theta \left[ \sum_{n=1}^{\xi} \beta^{n-1} X_n + \beta^\xi \frac{\lambda}{1-\beta} \right].$$

In [12], computational methods for calculating Gittins indices are described and applied to the normal, Bernoulli and negative exponential families with conjugate priors, which involve using backward induction to approximate the right-hand side of (2) with the supremum over $\xi \geq 0$ replaced by the supremum over $0 \leq \xi \leq N$ for some large horizon $N$. For $\beta$ close to 1, such computational methods become time consuming as a very large horizon $N$ is required to yield an accurate approximation. Thus it will be useful to have accurate analytic approximations to Gittins indices especially for $\beta$ close to 1.

In this paper, we consider the normal case with unknown mean $\theta$ and known variance where $\theta$ has a normal (conjugate) prior. Section 2 presents some monotonicity properties of the Gittins index. In particular, it is shown that the Gittins index is a nondecreasing function of the prior variance. In Section 3, a corrected diffusion approximation to the Gittins index is derived by embedding the (discrete-time) normal setting into the continuous-time Wiener process setting in which the Gittins index is determined by the stopping boundary for an optimal stopping problem (first introduced in [2]). By an application of Chernoff's continuity correction, the approximation includes a correction term which accounts for the difference between the discrete and continuous-time stopping boundaries. Numerical results are also given to assess the performance of this simple approximation. To prepare for the derivations, Sections 3.1 and 3.2 briefly review, respectively, some properties of the Gittins index for a Wiener process and Chernoff's continuity correction in optimal stopping.

The monograph of Gittins [12] provides a comprehensive theory of dynamic allocation indices and explores the class of problems whose optimal solutions can be characterized by dynamic allocation indices. On the other hand, Lai [17] and Chang and Lai [6] have proposed simple index-type adaptive allocation rules that are asymptotically optimal in both the Bayes and frequentist senses either as $N \to \infty$ (under uniform discounting) or as $\beta \to 1$ (under geometric discounting). Brezzi and Lai [5] have recently refined and modified these adaptive allocation rules in the presence of switching costs, while Hu and Wei [15] have constructed asymptotically optimal adaptive allocation rules subject to the irreversibility constraint. Various applications of the theory of multi-armed bandits can be found in sequential clinical trials, market pricing, labor markets and search problems; see e.g. [1, 8, 16, 19, 20].

## 2. Some monotonicity properties of the Gittins index for a normal reward process

In this section, we consider the case that $X_1, X_2, \ldots$ are (conditionally) iid $N(\theta, \sigma^2)$, the unknown mean $\theta$ has a prior $\pi = N(u, v)$ and the variance $\sigma^2$ is known. The Gittins index is denoted by $\lambda(u, v, \sigma^2, \beta)$. By location and scale equivariance properties (cf. [12], Section 6.4),

$$(3) \qquad \lambda(u, v, \sigma^2, \beta) = u + r\,\lambda(0,\, v/r^2,\, \sigma^2/r^2,\, \beta)$$

for $r > 0$.

**Lemma 1.** *The Gittins index $\lambda(u, v, \sigma^2, \beta)$ is nonincreasing in $\sigma^2$.*

*Proof.* We prove the lemma by a simple randomization argument. Fix $0 < \sigma_1^2 < \sigma_2^2$. Let $X_1, X_2, \ldots$ be (conditionally) iid $N(\theta, \sigma_1^2)$ given $\theta$, which is assumed to have a prior $\pi = N(u, v)$. Let $\epsilon_1, \epsilon_2, \ldots$ be iid $N(0, \sigma_2^2 - \sigma_1^2)$ (independent of the $X_i$). Then $X_1' = X_1 + \epsilon_1$, $X_2' = X_2 + \epsilon_2, \ldots$ are (conditionally) iid $N(\theta, \sigma_2^2)$ given $\theta$. For any stopping time $\xi' \geq 1$ with respect to the filtration $\mathcal{F}'$ generated by $X_1', X_2', \ldots$, we have

$$E_\pi E_\theta \sum_{n=1}^{\xi'} \beta^{n-1} X_n' = E_\pi E_\theta \sum_{n=1}^{\xi'} \beta^{n-1} X_n + E_\pi E_\theta \sum_{n=1}^{\infty} \beta^{n-1} \epsilon_n 1_{\{\xi' \geq n\}}$$

$$= E_\pi E_\theta \sum_{n=1}^{\xi'} \beta^{n-1} X_n.$$

Since every stopping time $\xi'$ with respect to $\mathcal{F}'$ may be viewed as a randomized stopping time with respect to $\mathcal{F}$ (the filtration generated by $X_1, X_2, \ldots$), it follows that

$$\lambda(u, v, \sigma_2^2, \beta) = \sup_{\xi' \geq 1} E_\pi E_\theta \left( \sum_{n=1}^{\xi'} \beta^{n-1} X_n' \right) \bigg/ E_\pi E_\theta \left( \sum_{n=1}^{\xi'} \beta^{n-1} \right)$$

$$\leq \sup_{\xi \geq 1} E_\pi E_\theta \left( \sum_{n=1}^{\xi} \beta^{n-1} X_n \right) \bigg/ E_\pi E_\theta \left( \sum_{n=1}^{\xi} \beta^{n-1} \right)$$

$$= \lambda(u, v, \sigma_1^2, \beta),$$

completing the proof. □

**Theorem 1.** $\lambda(0, v, \sigma^2, \beta)/\sqrt{v}$ *is nondecreasing in $v$.*

*Proof.* For fixed $0 < v_2 < v_1$, it follows from (3) and Lemma 1 that

$$\lambda(0, v_2, \sigma^2, \beta) = \sqrt{v_2/v_1}\, \lambda(0, v_1, \sigma^2 v_1/v_2, \beta)$$

$$\leq \sqrt{v_2/v_1}\, \lambda(0, v_1, \sigma^2, \beta),$$

completing the proof. □

**Corollary 1.** $\lambda(u, v, \sigma^2, \beta) = u + \lambda(0, v, \sigma^2, \beta)$ *is nondecreasing in $u$ and $v$.*

**Remark 1.** For the Wiener process setting, Bather [2] proved a result analogous to Theorem 1 (see (7) and (8) below).

**Remark 2.** For a normal two-armed bandit in which the means of arms 1 and 2 have independent normal priors $N(u_1, v_1)$ and $N(u_2, v_2)$ and their variances are known and equal, it follows from Corollary 1 that under geometric discounting, it is optimal to pull arm 1 initially if $u_1 \geq u_2$ and $v_1 \geq v_2$. It seems natural to conjecture that the same also holds under uniform discounting. Note that Berry [3] made a similar conjecture regarding a Bernoulli two-armed bandit, which has not been resolved (cf. [4], Section 7.3).

**Remark 3.** Along the lines of the proof of Theorem 1, it can be readily shown that

$$\lambda(0, v_1, \sigma_1^2, \beta)/\sqrt{v_1} \geq \lambda(0, v_2, \sigma_2^2, \beta)/\sqrt{v_2}$$

if $v_1 \geq v_2$ and $v_1/\sigma_1^2 \geq v_2/\sigma_2^2$. Note that for a normal distribution $N(\theta, \sigma^2)$ where $\theta$ has a normal prior $N(0, v)$, $v/\sigma^2$ may be referred to as the signal-to-noise ratio since $v$ is the second moment of the "signal" $\theta$.

## 3. Corrected diffusion approximation to the Gittins index for a normal reward process

In Section 3.3, we derive an approximation to the Gittins index for a normal distribution whose mean is assumed to have a normal prior. To prepare for the derivations, we briefly review, in Sections 3.1 and 3.2, some properties of the Gittins index for a Wiener process and Chernoff's continuity correction in optimal stopping.

### *3.1. Properties of the Gittins index for a Wiener process*

Bather [2] showed that for a Wiener process $\{W(t), t \geq 0\}$ with drift coefficient $\theta$ which has a normal prior $N(u_0, v_0)$, the Gittins index $\lambda^*(u_0, v_0, c)$ can be determined by the solution to an optimal stopping problem (to be described below) where $c > 0$ denotes the discount rate in continuous time (see also [6] and Section 6.6 of [12]). Here $\lambda^*(u_0, v_0, c)$ is defined as the infimum of the set of solutions $\lambda$ of the equation (cf. (2) )

$$
\begin{aligned}
\lambda \int_0^\infty e^{-ct} dt &= \sup_{\tau \geq 0} E_\pi E_\theta \left[ \int_0^\tau e^{-ct} dW(t) + \lambda \int_\tau^\infty e^{-ct} dt \right] \\
&= \sup_{\tau \geq 0} E_\pi \left[ \int_0^\tau \theta e^{-ct} dt + \lambda \int_\tau^\infty e^{-ct} dt \right] \\
&= \sup_{\tau \geq 0} E_\pi \left[ \int_0^\tau u(t) e^{-ct} dt + \lambda \int_\tau^\infty e^{-ct} dt \right] \\
&= \sup_{\tau \geq 0} E_\pi \left[ c^{-1} u_0 - c^{-1} \left( u(\tau) - \lambda \right) e^{-c\tau} \right],
\end{aligned}
$$
(4)

where the supremum is taken over all (real-valued) stopping times $\tau \geq 0$, $\pi = N(u_0, v_0)$ is the prior distribution of $\theta$, and $u(t)$ is the posterior mean of $\theta$, i.e.

(5) $$u(t) = E_\pi \left[ \theta \mid W(s), 0 \leq s \leq t \right] = \frac{v_0^{-1} u_0 + W(t)}{v_0^{-1} + t}.$$

The last equality in (4) follows from integration by parts along with the fact that a simple change of time transforms $u$ into standard Brownian motion, cf. $Y(v)$ below.

Define
$$v = v(t) = (v_0^{-1} + t)^{-1} \text{ (the posterior variance)}, \ s = v/c,$$
$$Y(v) = u_0 - u(t), \text{ and } Z(s) = Y(cs)/\sqrt{c}.$$

It can be readily shown that $\{Y(v), 0 < v \leq v_0\}$ is standard Brownian motion ( $Y(v_0) = 0$ ) in the $-v$ scale and $\{Z(s), 0 < s \leq s_0\}$ ($s_0 = v_0/c$) is standard Brownian motion ( $Z(s_0) = 0$ ) in the $-s$ scale. Letting $z_0 = (\lambda - u_0)/\sqrt{c}$, it follows that (4) is equivalent to

(6) $$z_0 e^{-1/s_0} = \sup_{0 < S \leq s_0} E \left[ \{Z(S) + z_0\} e^{-1/S} \right]$$

in the sense that $\lambda$ is a solution of (4) if and only if $z_0 = (\lambda - u_0)/\sqrt{c}$ is a solution of (6), where the supremum on the right-hand side of (6) is taken over all stopping times $0 < S \le s_0$ (in the $-s$ scale). It is more convenient to remove the restriction of $Z(s_0) = 0$ and rewrite (6) as

$$(6') \qquad z_0\, e^{-1/s_0} = \sup_{0 < S \le s_0} E\left[\, Z(S)\, e^{-1/S} \mid Z(s_0) = z_0 \,\right].$$

For the optimal stopping problem with payoff function $g(z,s) = ze^{-1/s}$ on the right-hand side of (6'), it is easily shown that the continuation region is of the form $\{(z,s) : z < b(s)\}$ where $b(s) > 0$ is the optimal stopping boundary. Since $z_0$ is a solution of (6') if and only if $(z_0, s_0)$ is in the stopping region (i.e. $z_0 \ge b(s_0)$), it follows that $\lambda^*(u_0, v_0, c)$, the infimum of the set of solutions $\lambda$ of the equation (4), satisfies $b(s_0) = \bigl(\lambda^*(u_0, v_0, c) - u_0\bigr)/\sqrt{c}$, i.e.

$$(7) \qquad \lambda^*(u_0, v_0, c) = u_0 + \sqrt{c}\, b(s_0) = u_0 + \sqrt{c}\, b(v_0/c).$$

Bather [2] showed that

$$(8) \qquad b(s)/\sqrt{s} \text{ is a nondecreasing function of } s,$$

$$(9) \qquad b(s) \le s/\sqrt{2} \text{ for all } s > 0, \text{ and } \lim_{s \to 0} b(s)/s = 1/\sqrt{2},$$

while Chang and Lai [6] derived the asymptotic expansion as $s \to \infty$

$$(10) \qquad b(s) = \left\{ 2s\left[\log s - \frac{1}{2}\log\log s - \frac{1}{2}\log 16\pi + o(1)\right] \right\}^{1/2}.$$

Based on (8)–(10) together with extensive numerical work (involving corrected Bernoulli random walk approximations for Brownian motion), Brezzi and Lai [5] have suggested the following closed-form approximation $\Psi(s)$ to $b(s)/\sqrt{s}$

$$(11) \qquad \frac{b(s)}{\sqrt{s}} \approx \Psi(s) = \begin{cases} \sqrt{s/2} & \text{for } s \le 0.2, \\ 0.49 - 0.11\, s^{-1/2} & \text{for } 0.2 < s \le 1, \\ 0.63 - 0.26\, s^{-1/2} & \text{for } 1 < s \le 5, \\ 0.77 - 0.58\, s^{-1/2} & \text{for } 5 < s \le 15, \\ \left\{ 2\log s - \log\log s - \log 16\pi \right\}^{1/2} & \text{for } s > 15. \end{cases}$$

### 3.2. Chernoff's continuity correction in optimal stopping

In his pioneering work, Chernoff [7] studied the relationship between the solutions of the discrete and continuous-time versions of the problem of testing sequentially the sign of the mean of a normal distribution. His result may be stated more generally as follows. Let $\{B(t)\}$ be standard Brownian motion and let $g(x,t)$ be a smooth payoff function for $t \le T$ (horizon) for which the continuation region is of the form $\{(x,t) : x < b(t)\}$. Consider a constrained optimal stopping problem where stopping is permitted only at $n\delta$, $n = 1, 2, \ldots$ where $\delta$ is a given (small) positive number. Suppose that there exist stopping boundary points $b_\delta(n\delta)$, $n = 1, 2, \ldots$ such that starting from $B(n_0\delta) = x_0$ for any given $n_0$ and $x_0$, the optimal stopping rule is to stop at the first $n \ge n_0$ at which $B(n\delta) \ge b_\delta(n\delta)$. So $b_\delta(n\delta)$ (or $b(t)$, resp.) is the discrete-time (or continuous-time, resp.) stopping boundary for the constrained (or unconstrained, resp.) optimal stopping problem. Then for fixed $t < T$, we have

$$(12) \qquad b_\delta(t) = b(t) - \rho\sqrt{\delta} + o(\sqrt{\delta}) \quad \text{as } \delta \to 0,$$

where $b_\delta(t) = b_\delta([t/\delta]\delta)$, $\rho = ES_{\tau_+}^2/2ES_{\tau_+} \approx 0.583$, $\tau_+ = \inf\{n : S_n > 0\}$, $S_n = X_1 + \cdots + X_n$, and the $X_i$ are iid $N(0,1)$.

Chernoff [7] derived (12) by relating the original problem to an associated stopping problem in which there is a horizon at $t = 0$ and the payoff function is $g(x,t) = -t + x^2 \, 1_{\{x<0,\, t=0\}}$, $t \le 0$. For the associated stopping problem, stopping is permitted at $0, -1, -2, \ldots$, and there exist stopping boundary points $b_{-1} > b_{-2} > \cdots$ such that starting from $(x_0, n_0)$ with $n_0 < 0$, the optimal stopping rule is to stop at the first $n_0 \le n \le 0$ at which

$$x_0 + X_1 + \cdots + X_{n-n_0} \ge b_n \quad (b_0 = -\infty).$$

Chernoff [7] and subsequently Chernoff and Petkau [9] and Hogan [14] showed that

$$\lim_{n \to -\infty} b_n = -ES_{\tau_+}^2 / 2\, ES_{\tau_+}$$

for normal, Bernoulli and general $X$ (with finite fourth moment), respectively. Recently, under mild growth conditions on $g$, Lai, Yao and AitSahlia [18] have proved (12) when the Brownian motion process is replaced by a general random walk in the constrained problem.

### 3.3. Approximating the Gittins index for a normal reward process

In this subsection, we consider the case that $X_1, X_2, \ldots$ are (conditionally) iid $N(\theta, \sigma^2)$ and the unknown mean $\theta$ has a prior $\pi = N(u_0, v_0)$. Without loss of generality, we assume $\sigma^2 = 1$. For notational simplicity, the Gittins index $\lambda(u_0, v_0, 1, \beta)$ will be abbreviated to $\lambda(u_0, v_0, \beta)$. Recall that $\lambda(u_0, v_0, \beta)$ is the infimum of the set of solutions $\lambda$ of the equation (2). As in Section 3.1, let $\{W(t), t \ge 0\}$ be a Wiener process with drift coefficient $\theta$ which has a normal prior $N(u_0, v_0)$. Noting that $(X_1, X_2, \ldots)$ and $(W(1), W(2) - W(1), \ldots)$ have the same joint distribution, we can rewrite (2) as

$$\begin{aligned}
\frac{\lambda}{1-\beta} &= \sup_{\xi \ge 0} E_\pi E_\theta \left[ \sum_{n=1}^{\xi} \beta^{n-1} \bigl( W(n) - W(n-1) \bigr) + \beta^\xi \frac{\lambda}{1-\beta} \right] \\
&= \sup_{\xi \ge 0} E_\pi \left[ \sum_{n=1}^{\xi} \beta^{n-1} u(n-1) + \beta^\xi \frac{\lambda}{1-\beta} \right] \\
&= \sup_{\xi \ge 0} E_\pi \left[ \frac{c}{1-\beta} \int_0^\xi u(t)\, e^{-ct} dt + \beta^\xi \frac{\lambda}{1-\beta} \right] \\
&= \frac{1}{1-\beta} \sup_{\xi \ge 0} E_\pi \left[ u_0 - \bigl( u(\xi) - \lambda \bigr) e^{-c\xi} \right],
\end{aligned}$$

where $u(t)$ is given in (5), $c = -\log \beta$ and the third equality follows since

$$\begin{aligned}
E &\left[ \frac{c}{1-\beta} \int_{n-1}^n u(t)\, e^{-ct} 1_{\{\xi \ge n\}}\, dt \, \Big| \, W(s),\, 0 \le s \le n-1 \right] \\
&= \frac{c}{1-\beta} 1_{\{\xi \ge n\}} \int_{n-1}^n u(n-1)\, e^{-ct}\, dt = \beta^{n-1} u(n-1) 1_{\{\xi \ge n\}}.
\end{aligned}$$

With the notation introduced in Section 3.1, we can further rewrite (2) as

$$\lambda - u_0 = \sup_{V \in \{v_0/(1+v_0 n),\, n=0,1,\ldots\}} E\left[\left(\lambda - u_0 + Y(V)\right) e^{-cV^{-1}+c v_0^{-1}}\right]$$

where the supremum is taken over all stopping times $V$ taking values in $\{v_0/(1+v_0 n),\, n = 0, 1, \ldots\}$. In terms of Brownian motion $Z(s)$ in the $-s$ scale, (2) is equivalent to

$$(13) \quad z_0 e^{-1/s_0} = \sup_{S \in \{c^{-1} v_0/(1+v_0 n),\, n=0,1,\ldots\}} E\left[Z(S)\, e^{-1/S} \mid Z(s_0) = z_0\right]$$

where $z_0 = (\lambda - u_0)/\sqrt{c}$, $s_0 = v_0/c$ and the supremum is taken over all stopping times $S$ taking values in $\{c^{-1} v_0/(1+v_0 n),\, n = 0, 1, \ldots\}$.

For the constrained optimal stopping problem on the right-hand side of (13), there exist optimal stopping boundary points $b_{v_0}(c^{-1} v_0/(1+v_0 n))$, $n = 0, 1, \ldots$ such that the optimal stopping rule is to stop at the first $n$ at which $Z(c^{-1}v_0/(1+v_0 n)) \geq b_{v_0}(c^{-1}v_0/(1+v_0 n))$. So $b_{v_0}(v_0/c)$ is the infimum of the set of solutions $z_0$ of the equation (13). It then follows that

$$(14) \quad \lambda(u_0, v_0, \beta) = u_0 + \sqrt{c}\, b_{v_0}(v_0/c).$$

Since in the constrained optimal stopping problem the permissible stopping time points $c^{-1} v_0/(1+v_0 n)$, $n = 0, 1, 2, \ldots$ are not equally spaced, there is no rigorous justification for applying (12) to relate the discrete and continuous-time stopping boundaries $b_{v_0}(t)$ and $b(t)$ for the constrained and unconstrained problems. However, it can be argued heuristically that (12) applies when the spacing between many successive permissible stopping time points is approximately constant ( cf. bottom of page 47 in [10]). Thus we arrive at the following approximation

$$(15) \quad b_{v_0}(v_0/c) \approx b(v_0/c) - 0.583\sqrt{\delta}$$

where

$$(16) \quad \delta = \frac{c^{-1}v_0}{1+v_0 \cdot 0} - \frac{c^{-1}v_0}{1+v_0 \cdot 1} = \frac{c^{-1}v_0^2}{1+v_0},$$

provided that $v_0/c$ is bounded away from $0$ (the horizon of the optimal stopping problem) and $\delta \approx \frac{c^{-1}v_0}{1+v_0 n} - \frac{c^{-1}v_0}{1+v_0(n+1)}$ for many (small to moderate) $n$'s. That is, we expect the approximation (15) to be reasonably good if $v_0$ is small and $v_0/c$ is not too close to 0. It follows from (14), (15), (16) and (11) that

$$(17) \quad \begin{aligned} \lambda(u_0, v_0, \beta) &\approx u_0 + \sqrt{c}\, b(v_0/c) - 0.583\, v_0/\sqrt{1+v_0} \\ &\approx u_0 + \sqrt{v_0}\, \Psi(v_0/c) - 0.583\, v_0/\sqrt{1+v_0}. \end{aligned}$$

Note that the continuation region for the constrained problem must be contained in the continuation region for the unconstrained problem, so that $b_{v_0}(v_0/c) < b(v_0/c)$. Thus the uncorrected diffusion approximation $u_0 + \sqrt{c}\, b(v_0/c)$ overestimates $\lambda(u_0, v_0, \beta) = u_0 + \sqrt{c}\, b_{v_0}(v_0/c)$, which is recorded in the following theorem.

**Theorem 2.** $\lambda(u_0, v_0, \beta) < u_0 + \sqrt{c}\, b(v_0/c)$ where $c = \log \beta^{-1}$.

A related upper bound for $\lambda(u_0, v_0, \beta)$ is given in Theorem 6.28 of Gittins [12], which states, in our notation, that

(18) $$\lambda(u_0, v_0, \beta) < u_0 + \sqrt{1-\beta}\ b(v_0/(1-\beta)).$$

Since $b(s)/\sqrt{s}$ is nondecreasing in $s$ by (8) and since $c = \log \beta^{-1} > 1 - \beta$, we have

$$\sqrt{c}\ b(v_0/c) \leq \sqrt{1-\beta}\ b(v_0/(1-\beta)),$$

so that the upper bound given in Theorem 2 is sharper than (18).

In the approximation (15), the correction term $0.583\sqrt{\delta}$ with $\delta$ given in (16) appears to be a little too large since the spacing between successive permissible stopping time points $\frac{c^{-1}v_0}{1+v_0 n} - \frac{c^{-1}v_0}{1+v_0(n+1)}$ is strictly less than $\delta$ for $n \geq 1$. To compensate for this overcorrection, we propose (in view of (9)) to replace $b(v_0/c)$ by $(v_0/c)/\sqrt{2}$ in (15), resulting in the following simple approximation

(19) $$\lambda(u_0, v_0, \beta) \approx u_0 + v_0/\sqrt{2c} - 0.583\, v_0/\sqrt{1+v_0}.$$

Note that (19) agrees with (17) for $v_0/c \leq 0.2$ in view of (11).

In his Table 1, Gittins [12] tabulates $n(1-\beta)^{1/2}\lambda(0, n^{-1}, \beta)$ for various values of $n$ and $\beta$. Our Table 1 compares $n(1-\beta)^{1/2}\lambda(0, n^{-1}, \beta)$ with the corrected and uncorrected approximations (based on (17) and (19))

(CA) $\quad n(1-\beta)^{\frac{1}{2}}\left[\dfrac{1}{\sqrt{n}}\Psi\left(\dfrac{1}{nc}\right) - \dfrac{0.583\, n^{-1}}{\sqrt{1+n^{-1}}}\right]$

$\quad\quad\quad = (1-\beta)^{\frac{1}{2}}\left[\sqrt{n}\,\Psi\left(\dfrac{1}{nc}\right) - \dfrac{0.583}{\sqrt{1+n^{-1}}}\right],$

(UA) $\quad n(1-\beta)^{\frac{1}{2}}\dfrac{1}{\sqrt{n}}\Psi\left(\dfrac{1}{nc}\right) = (1-\beta)^{\frac{1}{2}}\sqrt{n}\,\Psi\left(\dfrac{1}{nc}\right),$

(CA′) $\quad (1-\beta)^{\frac{1}{2}}\left[\dfrac{1}{\sqrt{2c}} - \dfrac{0.583}{\sqrt{1+n^{-1}}}\right],$

(UA′) $\quad \sqrt{(1-\beta)/(2c)}.$

**Remark 4.** As explained earlier, the uncorrected approximations have positive bias due to overestimation. The corrected approximations are reasonably accurate for moderate to large $n$ and for large $\beta$. For moderate $n$, (CA) (or (CA′), resp.) tends to underestimate (or overestimate, resp.) $n(1-\beta)^{1/2}\lambda(0, n^{-1}, \beta)$. This observation naturally leads to approximating $n(1-\beta)^{1/2}\lambda(0, n^{-1}, \beta)$ by the average of (CA) and (CA′), which is also included in Table 1. Overall, $[(CA) + (CA′)]/2$ has the best performance, while (CA′) is better than (CA) except for small $n$ and large $\beta$.

**Remark 5.** Table 1 of Gittins [12] suggests that $n(1-\beta)^{1/2}\lambda(0, n^{-1}, \beta)$ is increasing in $n$. For $\beta = 0.5, 0.6, 0.7, 0.8, 0.9, 0.95$, Gittins has numerically estimated $\lim_{n\to\infty} n(1-\beta)^{1/2}\lambda(0, n^{-1}, \beta)$. These numbers are compared in Table 2 with the limits $(1-\beta)^{1/2}[(2c)^{-1/2} - 0.583]$ (or $(1-\beta)^{1/2}/(2c)^{1/2}$, resp.) obtained from the corrected approximations (CA) and (CA′) (or uncorrected approximations (UA) and (UA′), resp.) as $n \to \infty$. It should be noted that the heuristic justification for the corrected approximations requires $v_0/c = 1/(n\,c)$ not to be very close to 0.

Table 1
*Gittins indices and approximations*
($\beta = 0.5, 0.7, 0.9, 0.95, 0.99, 0.995$)

|  | \multicolumn{5}{c}{$n$} |
|---|---|---|---|---|---|
|  | 10 | 50 | 100 | 500 | 1000 |
| $\beta = 0.5$ | | | | | |
| $n(1-\beta)^{1/2}\lambda(0, n^{-1}, \beta)$ | 0.211 | 0.224 | 0.226 | 0.227 | 0.227 |
| $[(CA) + (CA')]/2$ | 0.208 | 0.192 | 0.190 | 0.189 | 0.189 |
| (CA) | 0.208 | 0.192 | 0.190 | 0.189 | 0.189 |
| (CA$'$) | 0.208 | 0.192 | 0.190 | 0.189 | 0.189 |
| (UA) | 0.601 | 0.601 | 0.601 | 0.601 | 0.601 |
| (UA$'$) | 0.601 | 0.601 | 0.601 | 0.601 | 0.601 |
| $\beta = 0.7$ | | | | | |
| $n(1-\beta)^{1/2}\lambda(0, n^{-1}, \beta)$ | 0.311 | 0.337 | 0.341 | 0.344 | 0.345 |
| $[(CA) + (CA')]/2$ | 0.264 | 0.332 | 0.331 | 0.329 | 0.329 |
| (CA) | 0.184 | 0.332 | 0.331 | 0.329 | 0.329 |
| (CA$'$) | 0.344 | 0.332 | 0.331 | 0.329 | 0.329 |
| (UA) | 0.489 | 0.648 | 0.648 | 0.648 | 0.648 |
| (UA$'$) | 0.648 | 0.648 | 0.648 | 0.648 | 0.648 |
| $\beta = 0.9$ | | | | | |
| $n(1-\beta)^{1/2}\lambda(0, n^{-1}, \beta)$ | 0.415 | 0.480 | 0.493 | 0.504 | 0.506 |
| $[(CA) + (CA')]/2$ | 0.357 | 0.506 | 0.505 | 0.505 | 0.505 |
| (CA) | 0.201 | 0.506 | 0.505 | 0.505 | 0.505 |
| (CA$'$) | 0.513 | 0.506 | 0.505 | 0.505 | 0.505 |
| (UA) | 0.377 | 0.689 | 0.689 | 0.689 | 0.689 |
| (UA$'$) | 0.689 | 0.689 | 0.689 | 0.689 | 0.689 |
| $\beta = 0.95$ | | | | | |
| $n(1-\beta)^{1/2}\lambda(0, n^{-1}, \beta)$ | 0.425 | 0.519 | 0.540 | 0.562 | 0.566 |
| $[(CA) + (CA')]/2$ | 0.382 | 0.468 | 0.568 | 0.568 | 0.568 |
| (CA) | 0.190 | 0.367 | 0.568 | 0.568 | 0.568 |
| (CA$'$) | 0.574 | 0.569 | 0.568 | 0.568 | 0.568 |
| (UA) | 0.314 | 0.496 | 0.698 | 0.698 | 0.698 |
| (UA$'$) | 0.698 | 0.698 | 0.698 | 0.698 | 0.698 |
| $\beta = 0.99$ | | | | | |
| $n(1-\beta)^{1/2}\lambda(0, n^{-1}, \beta)$ | 0.353 | 0.499 | 0.549 | 0.618 | 0.633 |
| $[(CA) + (CA')]/2$ | 0.390 | 0.453 | 0.485 | 0.647 | 0.647 |
| (CA) | 0.130 | 0.257 | 0.322 | 0.647 | 0.647 |
| (CA$'$) | 0.650 | 0.648 | 0.647 | 0.647 | 0.647 |
| (UA) | 0.185 | 0.315 | 0.380 | 0.705 | 0.705 |
| (UA$'$) | 0.705 | 0.705 | 0.705 | 0.705 | 0.705 |
| $\beta = 0.995$ | | | | | |
| $n(1-\beta)^{1/2}\lambda(0, n^{-1}, \beta)$ | 0.304 | 0.457 | 0.516 | 0.614 | 0.638 |
| $[(CA) + (CA')]/2$ | 0.424 | 0.437 | 0.470 | 0.562 | 0.665 |
| (CA) | 0.181 | 0.209 | 0.274 | 0.458 | 0.665 |
| (CA$'$) | 0.667 | 0.665 | 0.665 | 0.665 | 0.665 |
| (UA) | 0.221 | 0.250 | 0.315 | 0.499 | 0.706 |
| (UA$'$) | 0.706 | 0.706 | 0.706 | 0.706 | 0.706 |

**Remark 6.** Brezzi and Lai [5] have proposed a simple approximation to Gittins indices for general distributions which is justified by making use of the functional central limit theorem as $\beta \to 1$. For Bernoulli distributions (with beta conjugate priors), their approximation provides fairly accurate results. When applied to normal distributions, their approximation reduces to the uncorrected approximation (UA). It will be of great interest to see whether and how Chernoff's continuity correction can apply to approximate Gittins indices for nonnormal distributions.

Table 2
*The limits of Gittins indices and approximations*

| $\beta$ | $\lim_{n\to\infty} n(1-\beta)^{1/2}\,\lambda(0, n^{-1}, \beta)$ | (CA) and (CA') | (UA) and (UA') |
|---|---|---|---|
| 0.5 | 0.227 | 0.189 | 0.601 |
| 0.6 | 0.283 | 0.257 | 0.626 |
| 0.7 | 0.345 | 0.329 | 0.648 |
| 0.8 | 0.417 | 0.409 | 0.669 |
| 0.9 | 0.509 | 0.505 | 0.689 |
| 0.95 | 0.583 | 0.568 | 0.698 |